\documentclass{article} 
\usepackage{amsmath,epsf}
\usepackage[latin1]{inputenc}
\usepackage{amsfonts}

\title{On certain spaces of \\lattice diagram polynomials}
\author{Jean-Christophe Aval\\ Laboratoire A2X, Universit\'e Bordeaux 1\\ 351 cours de la Lib\'eration, F-33405 Talence cedex\\ e-mail : {\tt aval@math.u-bordeaux.fr}}

\numberwithin{equation}{section}
\newtheorem{theorem}{Theorem}[section]
\newtheorem{proposition}{Proposition}[section]

\newtheorem{lemma}{Lemma}[section]

\newtheorem{definition}{Definition}[section]

\newtheorem{proof}{Proof}

\font \sc=cmr9

\def\C{{\mathbb C}}
\def\N{{\mathbb N}}

\def\Q{{\mathbb Q}}
\def\S{{\mathcal S}}
\def\T{{\mathcal T}}
\def\I{{\mathcal I}}
\def\F{{\mathcal F}}

\def\Young#1{\vbox{\smallskip\offinterlineskip
    \halign{&\vbox{##}\kern-\Thickness\cr #1}}}

\newdimen\Squaresize \Squaresize=20pt
\newdimen\Thickness \Thickness=.1pt
\newdimen\Correction \Correction=7pt

\def\Vide#1{\hbox{
       \vbox to \Squaresize{\vss
          \hbox to \Squaresize{\hss#1 \hss}\vss}
    \hskip-\Correction}
   \kern-\Thickness}

\def\Carre#1{\hbox{\vrule width \Thickness
   \vbox to \Squaresize{\hrule height \Thickness\vss
      \hbox to \Squaresize{\hss$\scriptstyle#1$\hss}
   \vss\hrule height\Thickness}
   \unskip\vrule width \Thickness}
   \kern-\Thickness}

\begin{document}

\maketitle

\begin{abstract}
The aim of this work is to study some lattice diagram polynomials $\Delta_D(X,Y)$ as defined in \cite{lattice} and to extend results of \cite{untrou}. We recall that $M_D$ denotes the space of all partial derivatives of $\Delta_D$. In this paper, we want to study the space $M^k_{i,j}(X,Y)$ which is the sum of $M_D$ spaces where the lattice diagrams $D$ are obtained by removing $k$ cells from a given partition, these cells being in the ``shadow'' of a given cell $(i,j)$ of the Ferrers diagram. We obtain an upper bound for the dimension of the resulting space $M^k_{i,j}(X,Y)$, that we conjecture to be optimal. These upper bounds allow us to construct explicit bases for the subspace $M^k_{i,j}(X)$ consisting of elements of $0$ $Y$-degree. 
\end{abstract}

\section{Introduction}

\begin{definition}
A {\sl lattice diagram} is a finite subset of $\N\times\N$.
For $\mu_1\geq \mu_2\geq \cdots \geq \mu_k>0$, we say that $\mu=(\mu_1,
\mu_2, \ldots,
\mu_k)$ is a {\sl partition} of $n$ if $n=\mu_1+\cdots +\mu_k$. We associate to
a partition $\mu$ its Ferrers diagram $\{(i,j)\,:\,0\le i\le k-1,\,
0\le j\le \mu_{i+1}\}$ and we use the symbol
$\mu$ for both the partition and its Ferrers diagram.
\end{definition}
Most definitions and conventions we use are similar to \cite{lattice}.
For example, given the partition $(4,2,1)$, its partition diagram is
  $$\Young{\Carre{2,0}\cr
           \Carre{1,0}&\Carre{1,1}\cr
           \Carre{0,0}&\Carre{0,1}&\Carre{0,2}&\Carre{0,3}\cr
           }\quad.$$
It consists of the lattice cells
$\{(0,0),(0,1),(0,2),(0,3),(1,0),(1,1),(2,0)\}$.

\begin{definition}
Given a lattice diagram $D=\{(p_1,q_1), (p_2,q_2),\ldots , (p_n,q_n)\}$
we define the {\sl lattice determinant}
  $$
  \Delta_D(X;Y)= \det \big\| x_i^{p_j}y_i^{q_j}\big\|_{i,j=1}^n\,,
  $$
where $X=X_n=\{x_1,x_2,\ldots,x_n\}$ and $Y=Y_n=\{y_1,y_2,\ldots,y_n\}$.
\end{definition}
The polynomial
$\Delta_D(X;Y)$ is  bihomogeneous of degree $|p|=p_1+\cdots +p_n$ in $X$
and of degree $|q|=q_1+\cdots +q_n$  in  $Y$. To insure that this definition
associates a unique polynomial to $D$ we require
that the list of lattice cells be given with respect to the following ``pseudo-lexicographic'' order: $\{(0,0),(1,0),(2,0),(0,1),(1,1),(0,2),(0,3)\}$.

For a polynomial $P(X;Y)$, the vector space spanned by all the partial
derivatives of $P$ of
all orders is denoted
${\mathcal L}_\partial[P]$. A permutation $\sigma\in \S_n$ acts diagonally on
a polynomial $P(X;Y)$ as follows:
$\sigma  P(X;Y)\,=\, P(x_{\sigma_1},x_{\sigma_2},\ldots
,x_{\sigma_n};y_{\sigma_1},y_{\sigma_2},\ldots ,y_{\sigma_n})$.
Under this action, $\Delta_D(X;Y)$ is clearly an alternant. It follows that for
any lattice diagram $D$ with $n$ cells, the vector space
$M_D = {\mathcal L}_\partial[\Delta_D(X;Y)]$
is an $\S_n$-module. Since $\Delta_D(X;Y)$ is  bihomogeneous, this module
affords a natural
bigrading. 

The most general problem exposed in \cite{lattice} and \cite{coltrouee} concerns the space $M_D$. The main question is to decide whether this space is $\S_n$-isomorphic to a sum of left regular representations or not. In the particular case where
$D$ corresponds to a partition $\mu$ the question leads to the ``$n!$ conjecture'' which asserts that the space
$M_\mu$ is a single
copy of the left regular representation. Many efforts to prove this conjecture were only sufficient to obtain it in some special cases (see \cite{allen}, \cite{aval}, \cite{gh}, \cite{orbit} for example).
In \cite{coltrouee}, the case where all the lattice cells of $D$ lies on a
single axis is solved.

The
next class of lattice diagrams that is of interest is obtained by removing a single cell from a partition diagram. Its interest comes in part from the fact that it gives a possible recursive approach for the $n!$ conjecture.
If $\mu$ is a partition of $n+1$, we  denote by $\mu/ij$ the lattice
diagram obtained by
removing the cell $(i,j)$ from the Ferrers diagram of $\mu$. We  refer to the cell
$(i,j)$ as
the {\sl hole} of  $\mu/ij$. It is conjectured in \cite{lattice} that the number of copies of
the left regular
representations in $M_{\mu/ij}$ is equal to the cardinality (which we denote by $s$) of the
$(i,j)$-{\sl shadow}, that is
the cardinality of
$\{(i',j')\in\mu : i'\ge i,\,j'\ge j\}$. 

A study of the subspace $M_{\mu/ij}(X)$ of $M_{\mu/ij}$ consisting of elements of 0 $Y$-degree can be found in \cite{untrou}, in which the corresponding ``four term recursion'' is proved by using the construction of explicit bases.

Here we study the following problem. Let $\mu$ be a partition of $n+k$. This partition is fixed and does not appear in the following notations. 

\begin{definition}
Let $M_{i,j}^k$ denote the following sum of vector spaces
$$M_{i,j}^k(X,Y)=\sum_{(a_1,b_1),\ldots,(a_k,b_k)}M_{\mu/\{(a_1,b_1),\ldots,(a_k,b_k)\}},$$
where the sum is over all the $k$-tuples of cells in the shadow of $(i,j)$. 
\end{definition}
Because of the ``shift'' operators (see \cite{lattice}, Proposition I.3 or section 2 in this paper) we have $M_{\mu/ij}=M_{i,j}^1$. 
Hence this space $M_{i,j}^k$ is a possible generalization of $M_{\mu/ij}$ if we want to make $k$ holes in the Ferrers diagram.
That is precisely this space that we want to study.
In fact we obtain an upper bound for the dimension of this object, that we conjecture to be optimal.

In the second section we introduce some ``shift'' operators which are useful to move the holes and the cells in the diagrams. The third section is devoted to the proof of an upper bound for the dimension of $M_{i,j}^k$ that is conjectured to be optimal. In the fourth section we study $M_{i,j}^k(X)$, the subspace of $M_{i,j}^k(X,Y)$ consisting of elements of 0 $Y$-degree, for which we obtain explicit bases.

\section{The ``shift'' operators}

For the sake of simplicity, we only settle the following propositions for $X$-shifts. Of course similar results also hold for $Y$-shifts.

\begin{proposition}
Let $L$ be a lattice diagram. Then for any integer $k\ge 1$ we have
$$p_k(\partial X)\Delta_L(X,Y)=\sum_{i=1}^n\epsilon(p_k(i;L))\Delta_{p_k(i,L)}(X,Y)$$
where ${p_k(i,L)}$ is obtained by replacing the $i$-th biexponent $p_i,q_i$ by $p_i-k,q_i$ and the coefficient $\epsilon(p_k(i;L))$ is different from zero only if the resulting diagram consists of $n$ distinct cells in the positive quadrant. Its sign is the sign of the permutation that reorders the obtained biexponents in increasing pseudo-lexicographic order.
\end{proposition}

\begin{proof}
This is a particular case of Proposition I.1 in \cite{lattice}.
\end{proof}

\noindent {\bf Remark 1.}
The diagram ${p_k(i,L)}$ is the diagram obtained by pushing down the $i$-th cell of $L$: its biexponent $(p_i,q_i)$ is replaced by $(p_i-k,q_i)$ which corresponds to $k$ steps down. The other biexponents are unchanged. This duality between the substractions on the set of biexponents and movements of cells of the diagram will be extensively employed all along this article, explicitly or implicitly. 

\begin{proposition}
Let $L$ be a lattice diagram. Then for any integer $k\ge 1$ we have
$$e_k(\partial X)\Delta_L(X,Y)=\sum_{1\le i_1<i_2<\cdots<i_k\le n}\epsilon(e_k(i_1,\ldots,i_k;L))\Delta_{e_k(i_1,\ldots,i_k;L)}(X,Y)$$
where $e_k(i_1,\ldots,i_k;L)$ is obtained by replacing the biexponents $(p_{i_1},q_{i_1}),\ldots,(p_{i_k},q_{i_k})$ by $(p_{i_1}-1,q_{i_1}),\ldots,(p_{i_k}-1,q_{i_k})$ and where the coefficient $\epsilon(e_k(i_1,\ldots,i_k;L))$ is a nonnegative integer which is different from zero only when the resulting diagram consists of $n$ distinct cells in the positive quadrant. 
\end{proposition}

\begin{proof}
The proof is almost the same as for the previous proposition. We write
$$e_k(X)=\sum_{1\le j_1<\cdots< j_k\le n}x_{j_1}\ldots x_{j_k}.$$
We develop the determinantal form of $\Delta_L$ with respect to the columns $j_1,\ldots,j_k$ to obtain the following expression where $\Delta_L^{i_1,\ldots,i_k}$ denotes the lattice diagram polynomial relative to the biexponents $i_1,\ldots,i_k$ of $L$ and $A_{i_1,\ldots,i_k;j_1,\ldots,j_k}$ the cofactor:
$$\Delta_L(X)=\sum_{1\le i_1<\cdots< i_k\le n}\Delta_L^{i_1,\ldots,i_k}(x_{j_1},\ldots ,x_{j_k}) A_{i_1,\ldots,i_k;j_1,\ldots,j_k}.$$
Next we derive to obtain 
\begin{align*}
\partial (x_{j_1}\ldots x_{j_k})\Delta_L(X)=\sum_{1\le i_1<\cdots< i_k\le n}&(c_{i_1,\ldots,i_k;j_1,\ldots,j_k}\Delta_{e_k(i_1,\ldots,i_k;L)}^{i_1,\ldots,i_k}(x_{j_1},\ldots ,x_{j_k}) \\
&\times A_{i_1,\ldots,i_k;j_1,\ldots,j_k}),
\end{align*}
where $c_{i_1,\ldots,i_k;j_1,\ldots,j_k}$ is a nonnegative integer. In fact $c_{i_1,\ldots,i_k;j_1,\ldots,j_k}$ appears to be independent of $j_1,\ldots,j_k$; therefore we can omit the subscript $j_1,\ldots,j_k$.
Thus we get
\begin{align*}
e_k(\partial X)\Delta_L(X)=\sum_{1\le i_1<\cdots< i_k< n}\sum_{1< j_1\le\cdots< j_k\le n}&(c_{i_1,\ldots,i_k}\Delta_{e_k(i_1,\ldots,i_k;L)}^{i_1,\ldots,i_k}(x_{j_1},\ldots ,x_{j_k})\\
&\times A_{i_1,\ldots,i_k;j_1,\ldots,j_k}).
\end{align*}
By recognizing the developpement of $\Delta_{e_k(i_1,\ldots,i_k;L)}$, we finally obtain the expected formula. The sign of the coefficient $\epsilon(e_k(i_1,\ldots,i_k;L))$ is the sign of the permutation that reorders the obtained biexponents in increasing pseudo-lexicographic order; in fact with the choice of the order this permutation is always the identity.
\end{proof}

Let $D$ be a lattice diagram. Let us denote by $D'$ its complementary in the positive quadrant (it is an infinite subset). We take the convention that for any diagram $D$, $\Delta_{D'}=\Delta_D$ (the prime name $D'$ will denote the complementary of a lattice diagram $D$).

\begin{proposition}\label{h}
With the previous notations, for any lattice diagram $L$ and for any integer $k\ge 1$, we have
$$h_k(\partial X)\Delta_L(X,Y)=\sum_{i_1<i_2<\cdots<i_k}\epsilon(h_k(i_1,\ldots,i_k,L'))\Delta_{h_k(i_1,\ldots,i_k,L')}(X,Y)$$
where $h_k(i_1,\ldots,i_k,L')$ is obtained by replacing the biexponents $(p_{i_1},q_{i_1}),\ldots,(p_{i_k},q_{i_k})$ of $L'$ by $(p_{i_1}+1,q_{i_1}),\ldots,(p_{i_k}+1,q_{i_k})$ and the coefficient $\epsilon(h_k(i_1,\ldots,i_k,L'))$ is a nonnegative integer which is different from zero only when the resulting $L'$ consists of distinct cells in the positive quadrant. 
\end{proposition}

\begin{proof}
We shall prove this proposition by induction on $k$. If $k=1$, then $h_1=e_1$ and the result is true since moving down a cell is equivalent to moving up a hole. Assume the result is true up to $k-1$. Then we use the fact that $h_k=h_{k-1}e_1-h_{k-2}e_2+\cdots+(-1)^kh_1e_{k-1}+(-1)^{k+1}e_k$. Indeed $h_{k-1}$ pushes up $(k-1)$ different holes. And $e_1$ moves down one cell, that is pushes up a hole. But if this hole has already been pushed up by $h_{k-1}$, we can view this situation as $k-2$ holes which have made one ``step'' and one which has made two steps. But this situation is killed by the term $-h_{k-2}e_2$. Thus by looking successively at the terms in the previous formula for $h_k$, we get what we want.
\end{proof}

\noindent {\bf Remark 2.}
One efficient application of Proposition \ref{h} is to give immediate proofs of Propositions 1-2-3-4 of \cite{aval} (these propositions provide a Groebner basis of the annulator ideal of $\Delta{\mu}$ when $\mu$ is a hook. The previous proofs were recursive and intricate but the results now become simple applications of Proposition \ref{h}.

\noindent {\bf Remark 3.}
In the particular case of two holes, we obtain
$$e_{k-1}(X)p_{l}(Y)(\partial)\Delta_{\mu/\{(i,j),(i+1,j)\}}=c_1\Delta_{\mu/\{(i,j),(i+k,j+l)\}}+c_2\Delta_{\mu/\{(i+l,j),(i,j+k)\}},$$
$$e_{l-1}(Y)p_{k}(X)(\partial)\Delta_{\mu/\{(i,j),(i,j+1)\}}=c_3\Delta_{\mu/\{(i,j),(i+k,j+l)\}}+c_4\Delta_{\mu/\{(i+l,j),(i,j+k)\}}$$
where $c_1$, $c_2$, $c_3$, $c_4$ are rational constants such that $c_1$ and $c_2$ are of the same sign and $c_3$ and $c_4$ are of opposite sign. 

By Proposition \ref{h} we can move simultaneously the two holes.
This implies that for any couple of holes $(h_1,h_2)$ in the shadow of $(i,j)$ then $\Delta_{\mu/\{h_1,h_2\}}\in M_{\mu/\{(i,j),(i,j+1)\}}+M_{\mu/\{(i,j),(i+1,j)\}}$ thus $$M_{i,j}^2=M_{\mu/\{(i,j),(i,j+1)\}}+M_{\mu/\{(i,j),(i+1,j)\}}.$$

The question of the generalization of the previous result when $k\ge 3$ appears spontaneously. Is it sufficient to take only the diagrams such that the holes form a partition of origin $(i,j)$ ? The answer is negative. For example it is easy to check that when $\mu=(3,2)$
$$\Delta_{\mu/\{(0,0),(1,0),(0,2)\}}\not \in M_{\mu/\{(0,0),(1,0),(0,1)\}}+M_{\mu/\{(0,0),(0,1),(0,2)\}}.$$

\section{The upper bound}

We define the annulator ideal of a vector subspace $M$ of $\Q[X_n,Y_n]$ as the following ideal: $I_M=\{P\in \Q[X_n,Y_n]\ :\ \forall Q\in M,\ P(\partial) Q=0\}$. If $M={\mathcal L}_\partial[P]$ then we denote its annulator ideal simply by $I_P$.
In the case of $M^k_{i,j}$ we denote $I_{M^k_{i,j}}$ by $I^k_{i,j}$.

We recall the following important result (\cite{gh}, Proposition 1.1): $M=I_M^{\perp}$, where the scalar product is defined by $(P,Q)=L_0(P(\partial) Q)$ and where $L_0$ is the linear form that associates to a polynomial its term of degree 0.

\subsection{About ideals}

We want here to prove the following

\begin{proposition}\label{I=I}
$$I^k_{i,j}=\bigcap_{(a_1,b_1),\ldots,(a_k,b_k)}I_{\partial x_{n+1}^{a_1}\partial y_{n+1}^{b_1}\cdots\partial x_{n+k}^{a_k}\partial y_{n+k}^{b_k}\Delta_{\mu}} \cap \Q[X_n,Y_n] \stackrel {\rm def} {=} {\I},$$
where the sum is over the $k$-tuples of different cells in the shadow of $(i,j)$ that we assume to be ordered in lexicographic order. 
\end{proposition}

\begin{proof}
By expanding $\Delta_{\mu}$ with respect to the last $k$ columns, we obtain:
\begin{align*}
\Delta_{\mu}(X_{n+k},Y_{n+k})=\sum_{(a_1,b_1),\ldots,(a_k,b_k)}&\pm\Delta_{\{(a_1,b_1),\ldots,(a_k,b_k)\}}(x_{n+1},\ldots,x_{n+k},y_{n+1},\ldots,y_{n+k})\\
&\times\Delta_{\mu/\{(a_1,b_1),\ldots,(a_k,b_k)\}}(X_n,Y_n).
\end{align*}
Thus for example:
$$\partial (x_{n+1}^{a_1}y_{n+1}^{b_1}\cdots x_{n+k}^{a_k}y_{n+k}^{b_k})\Delta_{\mu}(X_{n+k},Y_{n+k})=c\Delta_{\mu/\{(a_1,b_1),\ldots,(a_k,b_k)\}}(X_n,Y_n)+C$$
where $c$ is a rational constant (different from 0) and $C$ a linear combination with coefficients in $\Q[x_{n+1},y_{n+1},\ldots,x_{n+k},y_{n+k}]$ of $\Delta_{\mu/\{(a'_1,b'_1),\ldots,(a'_k,b'_k)\}}(X_n,Y_n)$ with $(a'_j,b'_j)$ in the shadow of $(a_j,b_j)$.

Hence we get what we want because:
\begin{itemize}
\item $\I\subset I_{i,j}^k$ follows from what precedes by looking at the constant term in $\Q[x_{n+1},y_{n+1},\ldots,x_{n+k},y_{n+k}]$;
\item $I_{i,j}^k\subset \I$ follows directly from what precedes.
\end{itemize}
\end{proof}

\subsection{Orbits}

The reasonning is inspired from \cite{lattice}, Theorem 4.2.

We consider two families $\alpha=(\alpha_1,\ldots,\alpha_h)$ and $\beta=(\beta_1,\ldots,\beta_l)$. To any injective tableau $T$ of shape $\mu$ with entries $1,\ldots,n+k$, we associate a point $(a(T),b(T))$ in $\C^{2n}$ by the classical process, ie: $a_i(T)=\alpha_{r_i(T)}$, $b_i(T)=\beta_{c_i(T)}$ where $r_i(T)$ (resp. $c_i(T)$) is the number of the row (resp. column) of $T$ where the entry $i$ lies in $T$. We define $\rho$ as the orbit of $(a,b)$ when $T$ varies over the $(n+k)!$ injective tableaux of shape $\mu$.
We introduce $J_{\rho}$ the ideal of polynomials that are zero all over the orbit, and $I=$gr$J_{\rho}$ and $H=I^{\perp}$ (we recall that gr is an operator that associate to a polynomial its term of maximum degree). It is now a classical result (cf. \cite{gh}, Theorem 1.1) that $I\subset I_{\Delta_{\mu}}$.

We now look at another orbit in $\C^{2n}$. We consider the set of tableaux $T$of shape $\mu$ with $n$ entries and $k$ white cells such that the $k$ white cells are in the shadow of $(i,j)$ (we denote this set of tableaux by ${\T}_{i,j}^k$).
By the same process as described above, we define an orbit $\rho^k$ in $\C^{2n}$.
Since the cardinality of ${\T}_{i,j}^k$ is ${s \choose k} n!$ the orbit $\rho^k$ has this cardinality. 
We introduce $J_{\rho^k}$ the ideal of polynomials that are zero all over the orbit, and $I^k=$gr$J_{\rho^k}$ and $H^k=(I^k)^{\perp}$. We have of course $\dim H^k={s \choose k} n!$.

We want to prove that $M_{i,j}^k\subset H^k$ and by \cite{gh}, Proposition 1.1, it is equivalent to prove that $I^k\subset {\I}.$

\subsection{Inclusion}

We want here to obtain the next proposition:

\begin{proposition}

We have the inclusion:
$$I^k\subset {\I}.$$
\end{proposition}

\begin{proof}
Let $P$ be a polynomial in $J_{\rho^k}$. Let us consider
$$Q(X_{n+k},Y_{n+k})=P(X_n,Y_n)\prod_{i'=1}^{i}(x_{n+1}-\alpha_{i'})\cdots\prod_{i'=1}^{i}(x_{n+k}-\alpha_{i'})\prod_{j'=1}^{j}(y_{n+1}-\alpha_{j'})\cdots\prod_{j'=1}^{j}(y_{n+k}-\alpha_{j'}).$$

We want to check that this polynomial is an element of $J_{\rho}$. We take an element $(\alpha,\beta)$ of $\rho$. If its projection on $\C^{2n}$ is in $\rho^k$ then $Q(\alpha,\beta)=0$ because of $P$. If not it must have at least one entry between $n+1$ and $n+k$ in the first $i$ rows or the first $j$ columns and we have still $Q(\alpha,\beta)=0$.
Thus gr$(P)\in  I_{\partial x_{n+1}^{i}\partial y_{n+1}^{j}\cdots\partial x_{n+k}^{i}\partial y_{n+k}^{j}\Delta_{\mu}}$.

For any set of k cells $\{(a_1,b_1),\ldots,(a_k,b_k)\}$ in the shadow of $(i,j)$, we observe that $\forall h,\  1\le h\le k$, $a_h\ge i$ and $b_h\ge j$. 
Hence gr$(P)$ is in $\I$, which was to be proved.
\end{proof}

\subsection{Conclusion}

The main result is now a consequence of all what precedes:

\begin{theorem} \label{ub}
If $\mu$ is a partition of $n+k$ and $s$ the cardinal of the shadow of the cell $(i,j)$, then we have:
$$\dim M^k_{i,j}\le {s \choose k} n!.$$
\end{theorem}

\noindent
{\bf Remark 4.}
Recall the proof of Theorem 1.1 of \cite{gh}: the previous reasonning implies that if equality holds in Theorem \ref{ub}, then $M^k_{i,j}$ decomposes as ${s \choose k}$ times the left regular representation.

\vskip 0.2 cm

Numerical examples and the fact that the construction described in the previous subsection affords the ``good'' upper bound in the case of one set of variables (see the next section) support the following conjecture.

\vskip 0.2 cm
\noindent {\bf Conjecture 1.} (F. Bergeron)
{\it
With the notations of the previous theorem:
$$\dim M^k_{i,j} = {s \choose k} n!.$$
}

\noindent {\bf Remark 5.}
When $k=1$, this conjecture reduces to Conjecture I.2 of \cite{lattice} and when $s=k$ or $k=0$ to the $n!$ conjecture.

\section{Case of one set of variables}

The goal of this section is to obtain an explicit basis for $M_{i,j}^k(X)$, ie the subspace of $M_{i,j}^k(X,Y)$ consisting of elements of 0 $Y$-degree.

\subsection{Construction}

We recall that in \cite{aval} is constructed a basis for $M_{\mu}^0$ made of monomial derivatives $M_S(\partial)\Delta_{\mu}$ of $\Delta_{\mu}$, where the objects $S$ varies over a set $\S(\mu)$ which depends on $\mu$. The corresponding construction can be found in \cite{aval}. The cardinality of $\S(\mu)$ is equal to $n!/\mu!$ where $\mu!=\mu_1!\mu_2!\cdots\mu_k!$. This cardinality is the number of injective, row-increasing tableaux of shape $\mu$.

Then we choose in the Ferrers diagram $\mu$ $k$ cells which are simultaneously in the shadow of $(i,j)$ and on the right edges of $\mu$. We denote by $\F_{\mu}^k$ the set of these objects. In the next figure, the chosen cells are cells with a circle and in the cell $(i,j)$ appears a $+$ sign. In this example $n=142$ and $k=10$.
\vskip 0.3 cm

\centerline{
\epsffile{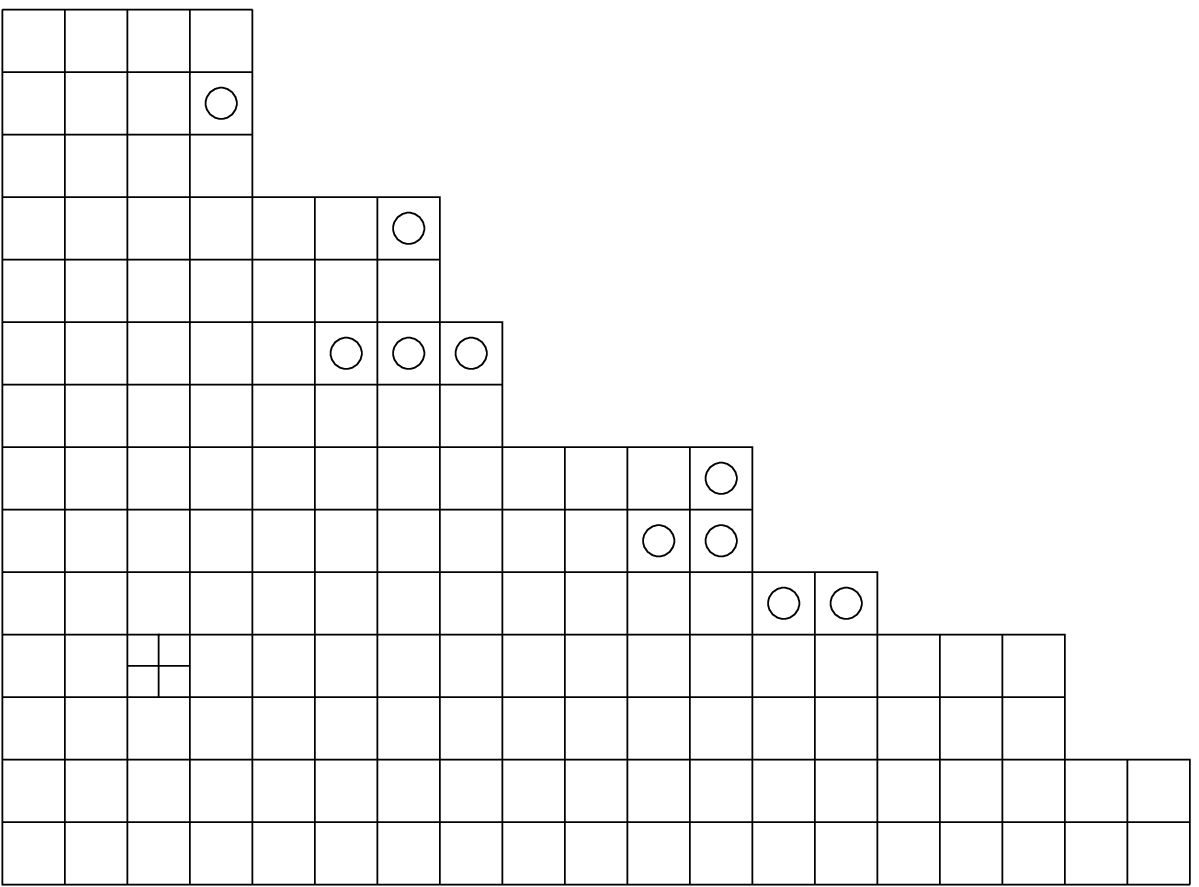}}

\vskip 0.3 cm

Then to an object $F$ in $\F_{\mu}^k$ we associate $\mu_{F}$ the partition of $n$ obtained by pushing up the circled cells and by removing the corresponding cells. We also define a partition $\mu_F^k$ with $k$ holes as follows. We look at the circled cells from the right to the left and from bottom to top. For a circled cell in column $j'\ge j$ such that there are $l$ possible places where this circle could have been placed in the same column and under this one, we put a hole in the cell $(i+l,j')$. The following figure illustrates this construction for the previous $F$ (the holes are as usual cells with crosses: $\times$).
\vskip 0.3 cm

\centerline{
\epsffile{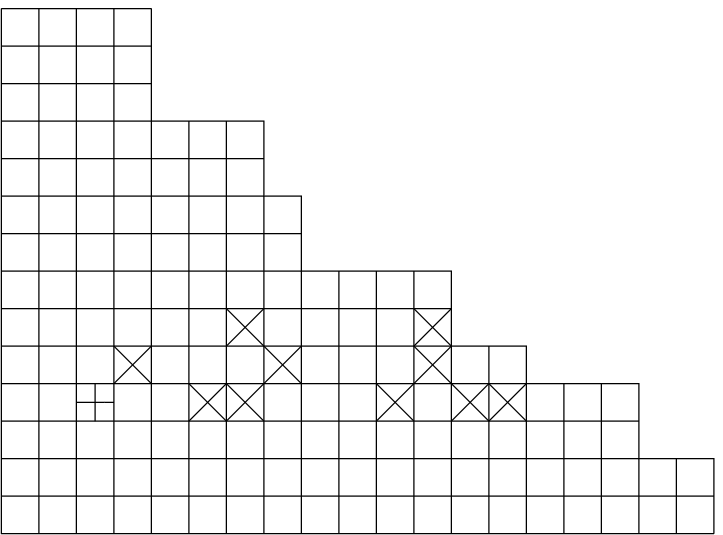}}

\vskip 0.3 cm

\begin{theorem}
With the previous notations
$$B_{i,j}^k(X)=\{M_S(\partial)\Delta_{\mu_F^k};\ S\in \S(\mu_F),\ F\in \F_{\mu}^k\}$$ is a basis for $M_{i,j}^k(X)$.
\end{theorem}

\subsection{Upper bound}

\begin{definition}

We denote by $\T_{i,j}^k$ the set of injective, row-increasing tableaux of shape $\mu$ with entries $1,\ldots,n$ (whence each of these tableaux has $k$ white cells) such that the $k$ white cells are in the shadow of $(i,j)$.  
\end{definition}

\begin{proposition}\label{upper bound}
We have the following upper bound for the dimension of $M_{i,j}^k(X)$:
$$\dim M_{i,j}^k(X)\le \#\T_{i,j}^k.$$
\end{proposition}

\begin{proof}

Let $I^k_{i,j}(X)$ denote the annulator ideal of the space $M^k_{i,j}(X)$. It is classical that by looking at the monomials in $Y_n$ on $\Q[X_n]$, we obtain $I^k_{i,j}(X)=I^k_{i,j}\cap\Q[X_n]$. From the Proposition \ref{I=I}, we deduce that $I^k_{i,j}(X)=I^k_{i,j}\cap\Q[X_n]=\I\cap\Q[X_n]$.

We look at the projection of the orbit $\rho^k$ on $\C^n$; it is equivalent to associate a point to each tableau in $\T_{i,j}^k$. Let $J_{\rho^k}^0$ denote the annulator ideal of this orbit.

The question is to prove the following inclusion: gr$(J_{\rho^k}^0)\subset I^k_{i,j}(X)$. 
The exactly same reasonning as in the case of two alphabets can be applied successfully.

\end{proof}

\subsection{Cardinality}

We claim that:
\begin{proposition}
We have the following equality
$$\#B_{i,j}^k(X)=\#\T_{i,j}^k.$$
\end{proposition}

\begin{proof}
It is sufficient to observe that we do not change the cardinality by pushing up the circled cells. We look at the example of the lines 9, 10 and 11 of the previous example.
\vskip 0.3 cm

\centerline{
\epsffile{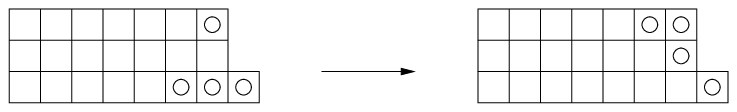}}

\vskip 0.3 cm
We observe that the lengths of the lines before the transformation are 5, 7, 6 and after the transformation 7, 6, 5. Thus the set of the lengths is unchanged. It is easy to see that it is always the case: the operation that pushes the holes up only permutes the lengths of the rows.
\end{proof}

\subsection{Independence}

We want here to conclude by proving that 

\begin{proposition}
The family defined in the paragraph Construction is linearly independent. Thus in particular equality holds in Proposition \ref{upper bound}.
\end{proposition}

\begin{proof}
Assume that we have a non-trivial dependence relation. 

We define the depth of a hole to be the number of cells (different from holes) that are above this hole.
We look at the $k$-tuples of the depths of the $k$ holes of the $\mu_F^k$: $(d_1\le d_2\le\ldots\le d_k)$. The crux of the proof is the following result:

\begin{lemma}
The $k$-tuples $(d_1,d_2,\ldots,d_k)$ are all distinct.
\end{lemma}

\begin{proof}
We want to prove that the depth of the holes increases from the right to the left and from top to bottom, and that two different $\mu_F^k$ give two different $k$-tuples of depths. We look at the circled cells with respect to this order.
We look at the next figure where $c$ denotes the number of circled cells in the column we consider, $m$ the number of positions under the lowest circled cell where we could put a circle (these cells appear with a square), $l$ the height of this column (we look only at the cells above the $i$-th row) and $l+h$ is the height of the ``next'' column (ie the first on the left).
 
\vskip 0.3 cm

\centerline{
\epsffile{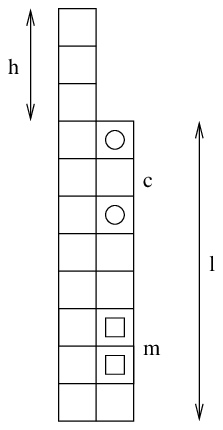}}

\vskip 0.3 cm

The depth of the lowest circled cell is $p=l-c-m$. The highest depth that could be obtained in this column is $l-c$ if $m=0$ and $l-c-1$ if $m>0$. And in the next column the lowest depth is (it corresponds to putting a circle at the top of the column): $l+h-1-c-h+1=l-c$. These observations are sufficient to conclude the proof of this lemma.
\end{proof}

Then if we have a non-trivial dependence relation, we consider the greatest $k$-tuple of depths with respect to the lexicographic order that appear in this relation (relative to an object $F^0$): $(d_1^0,d_2^0,\ldots,d_k^0)$. We then apply $h_k(\partial)^{d_1^0}.h_k(\partial)^{d_2^0-d_1^0}\ldots h_1(\partial)^{d_k^0-d_{k-1}^0}$.
It kills all the terms but those which come from the single object $F^0$. These terms give in fact terms in $B=\{M_S(\partial).\Delta_{\mu_{F^0}};\ S\in \S(\mu_{F^0})\}$, which are independent since $B$ is a basis of $M_{\mu_{F^0}}^0$.
\end{proof}

\noindent {\bf \large Acknowledgement.}
The author would like to thank François Bergeron for indicating him this problem and for numerous valuable suggestions.

\end{document}